\documentclass[a4paper,rotate,psfig,flegn]{article}
\usepackage[cp1251]{inputenc}
\usepackage [english]{babel}
\usepackage{amsmath}
\usepackage{amssymb}
\usepackage{graphicx}
\usepackage{subcaption}

\begin{document}
\bigskip
\centerline{\bf Spectral analysis of the Sturm-Liouville operator given on a system of segments}
\centerline{\bf S.Vovchuk}
\centerline{(V.N.Karazin Kharkiv National University)}
\medskip
\vskip5mm
\title{ Spectral analysis of the Sturm-Liouville operator given on a system of segments }


\title{ Spectral analysis of the Sturm-Liouville operator given on a system of segments }
\abstract {{}{} The spectral analysis of the Sturm-Liouville operator defined on a finite segment is the subject of an extensive literature  [1],[2]. Sturm-Liouville operators on a finite segment are well studied and have numerous applications [1,2]. The study of such operators already given on the system segments (graphs) was received in the works [3,4]. This work is devoted to the study of operators
$$(L_qy)(x)=col[-y_1''(x)+q_1(x)y_1(x), \ -y_2''(x)+q_2(x)y_2(x)],$$ where $y(x)=col[y_1(x),\ y_2(x)]\epsilon L^2(-a,0)\oplus L^2(0,b)=H, \ q_1(x), q_2(x) -$ real function $q_1\epsilon L^2(-a,0), q_2\epsilon L^2(0,b).$  Domain of definition $L_q$ has the form
$$\vartheta (L_q)={y=(y_1,y_2)\epsilon H; \ y_1\epsilon W_1^2(-a,0), \ y_2\epsilon W_2^2(0,b), \ y_1'(-a)=0, \ y_2'(b)=0; \ y_2(0)+py_1'(0)=0 \ y_1(0)+py_2'(0)=0}$$ $(p\epsilon \mathbb{R}, \ p\neq 0).$ Such an operator is self-adjoint in $H.$ The work uses the methods described in work  [5]. The main result is as follows: if the $q_1, q_2$ are small  (the degree of their smallness is determined by the parameters of the boundary conditions and the numbers $a, b$), ), then the eigenvalues  $\{\lambda_k(0)\}$ of the unperturbed operator $L_0$ are simple, and the eigenvalues  $\{\lambda_k(q)\}$ of the perturbed operator $L_q$ are also simple and located small in the vicinity of the points $\{\lambda_k(0)\}$.

  }

 \section{ Unperturbed operator }

  Consider the Hilbert space $H=L^2(-a,0)\oplus L^2(0,b)$, $(a,b\geq 0)$  by vector functions  $y(x)=col[y_1(x),\ y_2(x)],$ where $ y_1\epsilon L^2(-a,0), y_2\epsilon L^2(0,b).$ Define in $H$ a linear operator
\numberwithin{equation}{section}
\begin{equation}\label{Operator}\left(L_qy\right)(x)=col\left[
     -y_1''(x)+q_1(x)y_1(x), \\
    -y_2''(x)+q_2(x)y_2(x) \\
\right],
\end{equation}
where $q_1,\ q_2-$ real function and  $q_1(x)\epsilon L^2(-a,0), \ q_2\epsilon L^2(0,b). $

The domain of definition of the operator $L_q$ has the form,
\begin{equation} \label{Domain} \begin{split}\vartheta(L_q)= \{y=col[y1,y2]\epsilon H;\ y_1\epsilon W_1^2(-a,0),\
y_2\epsilon W_2^2(0,b),\ y_1'(-a)=0,\ y_2'(b)=0, \ y_2(0)+py_1(0)=0, \\ y_1'(0)+py_2'(0)=0\}
\end{split}
\end{equation}
 $ (p\epsilon\mathbb R, p\neq 0).$
 The operator $L_q$ (\ref{Operator}),(\ref{Domain}) is symmetric because
 \begin{eqnarray*}
 \left< L_qy,g\right>-\left<y,L_qg\right>=-\left.y'_1(x)\bar{g_1}(x)\right|_{-a}^0 +\left.y_1(x)\bar{g_1}'(x)\right|_{-a}^0-\left.y'_2(x)\bar{g_2}(x)\right|_0^b+\left.y_2(x)\bar{g_2}'(x)\right|_0^b=
 \\=py'_2(0)\bar{g_1}(0)-p\bar{g_2}'(0)y_1(0)-p\bar{g_1}(0)y'_2(0)+py_1(0)\bar{g_2}'(0)=0.
 \end{eqnarray*}
It is easy to show that  $L_q$ (\ref{Operator}),(\ref{Domain}) is self-adjoint.

Primary study the unperturbed operator  $L_0$ $(q_1=q_2=0)$. Respectively, the function of operator
$L_0$ is a solution to the equations
  \begin{equation} \label{Unperturbed} -y_1''=\lambda^2y_1, \  -y_2''=\lambda^2y_2 \quad (\lambda \epsilon \mathbb{C})
\end{equation}
and satisfies the boundary conditions  (\ref{Domain}). From the first two boundary conditions we find that
 \begin{equation} \label{Solutions} y_1=A\cos\lambda(x+a), \  y_2=B\cos\lambda(x-b),
\end{equation}
where $A, B,  C \epsilon \ \mathbb C.$
Second and third boundary  $\{y_1,y_2\}$  (\ref{Solutions}) give a system of equations for  $A $ and $ B$,
\begin{equation}\label{Sys}
 \begin{cases}
   Ap\cos\lambda a+B\cos\lambda b=0,\\
   -A\lambda\sin\lambda a+B p\lambda \sin\lambda b=0.
 \end{cases}
\end{equation}

This system has non-trivial solution $A$ and $B$, if only its determinant  $\bigtriangleup(0,\lambda)=0$,
  where
\begin{equation}\label{Det}
 \bigtriangleup(0,\lambda)=\lambda p^2\cos\lambda a\sin\lambda b+\lambda \cos\lambda b \sin\lambda a.
\end{equation}

If  $\lambda=0,$ then $y_1=A, \ y_2=B$ and from $y_2(0)+py_1(0)=0$ follows that $B=-pA.$  So $y=Acol[1,-p] \ (A\epsilon \mathbb{C})$ operator's own function  $L_0,$ responding to its own value  $\lambda=0.$
With $\lambda\neq 0$  from $\bigtriangleup(0,\lambda)=0$ follows
 \begin{equation} \label{Kh.equation} p^2\cos\lambda a\sin\lambda b+ \cos\lambda b \sin\lambda a=0,
\end{equation}
\newtheorem{remark}{Remark}
\begin{remark}
If $p=\pm1,$ that from  (\ref{Kh.equation}) follow $\sin\lambda(a+b)=0$
and hence the own numbers have the form \begin{equation}\label{Lambda}
 \lambda_n =\frac{\pi n}{a+b} \quad  (n\epsilon \mathbb{Z})
\end{equation}
\end{remark}

Consider the general case, not assuming, that  $p\pm1$ and write equality  (\ref{Kh.equation}) in form
$$(p^2+1)\sin\lambda(a+b)- (1-p^2)\sin\lambda(b-a)=0$$
or  \begin{equation}\label{Eq}
 \sin\lambda(a+b)-\frac{(1-p^2)}{(1+p^2)}\sin\lambda(b-a)=0
\end{equation}
\newcommand{\eqdef}{\stackrel{\mathrm{def}}{=}}
Let be  \begin{equation} \label{W} \lambda(b+a)=w, \ \frac{1-p^2}{1+p^2}=k, \ \frac{b-a}{b+a}=q,
\end{equation}
  then it is obvious that $|k|\leq 1, \ |q|\leq 1,$ and the equation  (\ref{Eq}) has form
 \begin{equation} \label{Sin} f(w)=0;\ f(w)\eqdef \sin w-k\sin q w
\end{equation}
The function $f(w)$ is odd, therefore it is enough to find its zeros $f(w)$  on the ray $\mathbb R_+.$

Show that the zeros of  $f(w)$ are simple. Assuming the opposite, suppose that $w-$ repeated root, then from $f(w)$ and $f'(w)=0,$
follows, that
$$   \begin{cases}
  \sin w = k \sin q w \\
   \cos w =k q\cos q w
 \end{cases}$$
 That means
$$k^2 \sin^2 q w+k^2q^2\cos^2 q w=1 $$
that's why
\begin{equation} \label{fif'}
 \sin^2 q w+q^2\cos^2 q w=\frac{1}{k^2}
\end{equation}

Since $|k|\leq 1 \ (k=1\ \text{{which}\ p=0, \text  is impossible by assumption})$ and $|q|< 1,$ then from  (\ref{fif'}) follows that the left side $\sin^2 q w+q^2\cos^2 q w\leq1,$ and right side  $\frac{1}{k^2}>1.$ That's why roots  $f(w)$ are simple.
 \newtheorem{Th}{Theorem }
\begin{Th}
 Roots $\{\lambda_s(0)\}$  of the characteristic function
 $\bigtriangleup(0,\lambda)$ (\ref{Det}) are simple except $\lambda_0(0)=0,$  which is duble multiple and they have the form,
\begin{equation}\label{Th}\Lambda_0=\{ 0,\ \lambda_s(0)=\pm\frac{w_s}{a+b},\ w_s>0;\ \sin w_s=k\sin q w_s \},
\end{equation}
where  $k \ i \ q-$ have of form  (\ref{W}), and the numbers  $w_s\epsilon \mathbb R_+ $ are numbered in ascending order.
 \end{Th}
\begin{remark}\label{remark-2}
Greatest positive root  $w_1$ equation $f(w)=0$ obviously lies in the interval  $\frac{\pi}{2}<w_1<\pi$ and that mean  $\frac{\pi}{2(a+b)}<\lambda_1(0)<\frac{\pi}{a+b}.$
\end{remark}

Eigenfunctions $\varphi(0,\lambda_s(0))$ of operator  $L_0,$ responding $\lambda_s(0)\epsilon\Lambda_0$ (\ref{Th})  are equal
\begin{equation}\varphi (0,\lambda_s(0))=A_s col\left[\cos\lambda_s(0) b \cos\lambda_s(0) (x+a), \ -p\cos\lambda_s(0) a \cos\lambda_s(0) (x-b)\right],
\end{equation}
which is an obvious consequence  (\ref{Solutions}), (\ref{Sys})

\section{Perturbed operator}

  Let's move on to the perturbed operator  $L_q.$ The equation for the eigenfunction  $y=col[y_1, y_2]$ of operator  $L_q$ has the form
\begin{equation}\label{L_q}
-y''_1+q_1y_1=\lambda^2 y_1, \\
-y''_2+q_2y_2=\lambda^2 y_2
\end{equation}

  Consider the integral equations
 \begin{equation}\label{Integral}
 \begin{cases}
y_1(x)=A\cos\lambda(x+a)+\int\limits_a^x \frac{\sin\lambda (x-t)}{\lambda}q_1(t)y_1(t)\,dt; \\
y_2(x)=B\cos\lambda(x-b)-\int\limits_x^b \frac{\sin\lambda (x-t)}{\lambda}q_2(t)y_2(t)\,dt.
\end{cases}
\end{equation}
 Then $\{y_k(x)\}-$ ) solution (\ref{Integral}) satisfy the equations (\ref{L_q}), and the first boundary conditions  (\ref{Domain}) correspond to  $y_1, y_2$.

Solvability of the integral equation  (\ref{Integral}) for $y_1.$  Definition of Volterra operator in  $L^2(-a,0),$
 \begin{equation}
 (K_1f)(x)=\int\limits_{-a}^x K_1(x,t)q_1(t)f(t)\,dt \quad (f\epsilon L^2(-a,0)),
 \end{equation}
 where
 \begin{equation}\label{K1f}
 (K_1f)(x)=\frac{\sin\lambda(x-t)}{\lambda}.
 \end{equation}
 Then the first of the equations in (\ref{Integral})  will take the form
 \begin{equation}
 (I-K_1)y_1=A\cos\lambda(x-t),
 \end{equation}
 And that means
 \begin{equation}\label{Y1}
 y_1=\sum_{n=0}^\infty K_1^nA\cos\lambda(x+a)
 \end{equation}
 where \begin{equation}
 (K_1^nf)(x)=\int\limits_{-a}^x K_{1,n}(x,t)q_1(t)f(t)\,dt ,
 \end{equation}
 For cores  $K_{1,n}(x,t)$ the recurrence relations are valid
  \begin{equation}\label{K1n}
 K_{1,n+1}(x,t)=\int\limits_t^x K_1(x,s)K_{1,n}(s,t)q_1(s)\,dt \quad (n> 1)
 \end{equation}
 where $K_1(x,t)$ have from  (\ref{K1f})
 \newtheorem{Lem1}{Lemma}
\begin{Lem1}
The kernels  $K_{1,n}(x,t)$ (\ref{K1n}) satisfy the inequalities
 \begin{equation}\label{Mark}
 |K_{1,n}(x,t)|\leq ch \beta(x-t)\frac{(x-t)^n}{n^n}\cdot\frac{\sigma_1^{n-1}(x)}{(n-1)!},
 \end{equation}
where  \begin{equation}
\beta=Im \lambda, \quad \sigma_1(x)=\int\limits_{-a}^x |q_1(t)|\,dt .
 \end{equation}
\end{Lem1}
The proof of the estimates  (\ref{Mark}) is carried out by induction
From (\ref{Y1}) it follows, that
$$y_1(\lambda,x)=A\cos\lambda(x+a)+A\int\limits_{-a}^x \sum_{n=1}^\infty K_{1,n}(x,t)q_1(t)\cos\lambda(t+a)\,dt=A\cos\lambda(x+a)+A\int\limits_{-a}^x N_1(x,t,\lambda)q_1(t)\cos\lambda(t+a)\,dt ,$$
where 
$$ N_1(x,t,\lambda)= \sum_{n=1}^\infty K_{1,n}(x,t).$$
 it follows from the estimates (\ref{Mark}) that this series converges and
 $$|N_1(x,t,\lambda)|\leq \cosh\beta(x-t)(x-t)\exp\left[(x-t)\sigma_1(x)
 \right]$$
Similar reasoning is valid for the second equation (\ref{Integral}).
\begin{Th}{}
Integral equations  (\ref{Integral}) are resolved and,-
\begin{equation}\label{IntSol}
 \begin{cases}
y_1(\lambda,x)=A\left(\cos\lambda(x+a)+\int\limits_{-a}^x N_1(x,t,\lambda)q_1(t)\cos\lambda(t+a)\,dt\right); \\
y_2(\lambda,x)=B\left(\cos\lambda(b-x)-\int\limits_{x}^b N_2(x,t,\lambda)q_2(t)\cos\lambda(b-t)\,dt\right),
\end{cases}
\end{equation}
In this case, the kernels  $\{N_k(x,t,\lambda)\}$ ) satisfy the estimates \begin{equation}\label{N_k}
 |N_k(x,t,\lambda)|\leq \cosh \beta(x-t)\cdot(x-t)\cdot exp\{(x-t)\sigma_k(t)\} \ (k=1,2),
 \end{equation}
where \begin{equation}
\beta=Im \lambda, \quad \sigma_1(x)=\int\limits_{-a}^x |q_1(t)|\,dt, \quad \sigma_2(x)=\int\limits_{x}^b |q_2(t)|\,dt .
 \end{equation}
\end{Th}

To find a characteristic function $\bigtriangleup(q,\lambda)$ the operator  $L_q$ L uses the last boundary conditions  (\ref{Domain}) for the  $\{y_k(\lambda, x)\},$ ), as a result, we obtain a one-row system of equations for  $A$ and $B$,-
\begin{equation}\label{Sys1}
 \begin{cases}
pA\left(\cos\lambda a+\int\limits_{-a}^0 N_1(0,t,\lambda)q_1(t)\cos\lambda(t+a)\,dt\right)+B\left(\cos\lambda b-\int\limits_{0}^b N_2(0,t,\lambda)q_2(t)\cos\lambda(b-t)\,dt\right)=0, \\
A\left(-\lambda \sin\lambda a+\int\limits_{-a}^0 N_1^\prime(0,t,\lambda)q_1(t)\cos\lambda(t+a)\,dt\right)+pB\left(\lambda\sin\lambda b-\int\limits_{0}^b N_2^\prime(0,t,\lambda)q_2(t)\cos\lambda(b-t)\,dt\right)=0
\end{cases}
\end{equation}
System (\ref{Sys1}) at this value  $q_1=q_2=0$ coincides with the system  (\ref{Sys}) and it has a nontrivial solution  $A, \ B,$ if its determinant
 $\bigtriangleup(q,\lambda)=0,$ where
\begin{equation}\label{Det1}
 \bigtriangleup(q,\lambda)\eqdef\begin{vmatrix} \ p\left(\cos\lambda a+\int\limits_{-a}^0 N_1(0,t,\lambda)q_1(t)\cos\lambda(t+a)\,dt\right) & \ \cos\lambda b-\int\limits_{0}^b N_2(0,t,\lambda)q_2(t)\cos\lambda(b-t)\,dt \\
-\lambda \sin\lambda a+\int\limits_{-a}^0 N_1^\prime(0,t,\lambda)q_1(t)\cos\lambda(t+a)\,dt &\ p\left(\lambda\sin\lambda b-\int\limits_{0}^b N_2^\prime(0,t,\lambda)q_2(t)\cos\lambda(b-t)\,dt\right) \end{vmatrix}
\end{equation}
It follows that,
\begin{equation}\label{delta}
\bigtriangleup(q,\lambda)=\bigtriangleup(0,\lambda)+\Phi (\lambda),
\end{equation}
where $\bigtriangleup(0,\lambda)$ ) have form  (\ref{Det}), and $\Phi (\lambda)$ is equal
\begin{equation}\label{Phi}
\begin{split}
\Phi (\lambda)=  p^2\{\lambda \sin \lambda b\int\limits_{-a}^0 N_1(0,t,\lambda)q_1(t)\cos\lambda(t+a)\,dt-\cos\lambda a\int\limits_{0}^b N_2^\prime(0,t,\lambda)q_2(t)\cos\lambda(b-t)\,dt -\\ -\int\limits_{-a}^0 N_1(0,t,\lambda)q_1(t)\cos\lambda(t+a)\,dt \times \int\limits_{0}^b N_2^\prime(0,t,\lambda)q_2(t)\cos\lambda(b-t)\,dt\} -\\-\lambda \sin\lambda a\int\limits_{0}^b N_2(0,t,\lambda)q_2(t)\cos\lambda(b-t)\,dt-\cos\lambda b\int\limits_{-a}^0 N_1^\prime(0,t,\lambda)q_1(t)\cos\lambda(t+a)\,dt+\\ +\int\limits_{0}^b N_2(0,t,\lambda)q_2(t)\cos\lambda(b-t)\,dt\cdot \int\limits_{-a}^0 N_1^\prime(0,t,\lambda)q_1(t)\cos\lambda(t+a)\,dt
\end{split}
\end{equation}
\begin{Th}
Operator characteristic function  $\bigtriangleup(q,\lambda)$ (\ref{Det1}) is expressed in terms of the operator  $L_q$ (\ref{Operator}),(\ref{Domain}) characteristic function  $\bigtriangleup(0,\lambda)$ (\ref{Det})  $L_0$ $(q_1=q_2=0)$ ) by the formula (\ref{delta}), where $\Phi (\lambda)$ has the form  (\ref{Phi}) and is an entire function of exponential type while it satisfies the estimate
\begin{equation}\label{phi}
|\Phi (\lambda)|\leq ch \beta a\cdot ch \beta b \cdot (\delta_1|\lambda|+\delta_2),
\end{equation}
 where
 \begin{equation}\label{2.19}\delta_1\eqdef \sigma_1 a e^{\sigma_1 a}+ \sigma_2 b e^{\sigma_2b}, \ \delta_2\eqdef \sigma_1 e^{\sigma_1 a}+ \sigma_2 e^{\sigma_2b}+\sigma_1 \sigma_2 (a+b)e^{\sigma_1 a+\sigma_2 b}
  \end{equation}
  and  $ \beta=Im \lambda, \  \sigma_1=\sigma_1(0),\  \sigma_2=\sigma_2(0)$.
 \end{Th}
 \newenvironment{Proof}
{\par\noindent{\bf Proof}}
{\hfill$\scriptstyle\blacksquare$}
  \begin{Proof}
 The estimates are similarly (\ref{N_k})  valid
 $$\left|\frac{\partial}{\partial x} N_k(x,t,\lambda)\right|\leq ch\beta(x-t)exp\{\sigma_k(x)(x-t)\} \ (k=1,2),$$
 therefore, it follows from  (\ref{Phi}) that
 \begin{eqnarray*}
 |\Phi(\lambda)|\leq p^2\{|\lambda|ch\beta b\cdot \cos\beta a\cdot e^{\sigma_1 a}\sigma_1 a+ch\beta a\cdot ch\beta b\cdot e^{\sigma_2 b}\sigma_2 +a\cdot ch\beta a\cdot ch\beta b\cdot e^{\sigma_1 a}\cdot e^{\sigma_2 b}\sigma_1 \sigma_2 \}+\\+|\lambda|ch\beta a ch\beta b e^{\sigma_2 b}\sigma_2 b+ ch\beta b\cdot ch\beta a\cdot e^{\sigma_1 a}\sigma_1+b ch\beta a\cdot ch\beta b\cdot e^{\sigma_1 a}\cdot e^{\sigma_2 b}\sigma_1\sigma_2.
 \end{eqnarray*}
  Thus,
  \begin{eqnarray*}|\Phi(\lambda)|\leq ch\beta b\cdot ch\beta a \{\sigma_1\cdot e^{\sigma_1 a}(1+|\lambda|p^2a)+\sigma_2\cdot e^{\sigma_2 b}(b|\lambda|+p^2)+\sigma_1 \sigma_2 e^{\sigma_1 a+\sigma_2 b}(b+p^2a)\}  \end{eqnarray*}
  And since  $p^2<1,$ ) then
  $$ |\Phi(\lambda)| \leq ch\beta a\cdot ch\beta b\{|\lambda|(\sigma_1 a e^{\sigma_1 a}+b\sigma_2 e^{\sigma_2 b})+\sigma_1 e^{\sigma_1 a}+\sigma_2 e^{\sigma_2 b}+\sigma_1 \sigma_2e^{\sigma_1 a+\sigma_2 b}(b+a)\}$$
  which proves  (\ref{phi}).
  \end{Proof}
 \section{Basic assessments}
  Characteristic function  $\bigtriangleup(0,\lambda)$ (\ref{Det}) taking into account these  (\ref{Kh.equation}),(\ref{Lambda}) is equal to,-
  \begin{equation}\label{3.1}
  \bigtriangleup(0,\lambda)=\lambda(p^2+1)Q(\lambda); \ Q(\lambda)\eqdef \sin\lambda (a+b)-k\sin q\lambda (a+b),
  \end{equation}
   where $q, k$ has form  (\ref{W}) and $|k|\leq1,\ |q|\leq1.$ Let us expand  $Q(\lambda)$ by the Taylor formula in a real neighborhood of the point $\lambda_s(0) (\neq0)$ (\ref{Th}),
  $$Q(\lambda)=(\lambda-\lambda_s)Q^\prime (\lambda_s)+\frac{(\lambda-\lambda_s)^2}{2}Q^{\prime\prime}(\xi_s)=(\lambda-\lambda_s)Q^\prime (\lambda_s)\left(1+\frac{(\lambda-\lambda_s)^2}{2}\cdot\frac{Q^{\prime\prime}(\xi_s)}{Q^\prime (\lambda_s)}\right),$$
  where $\lambda\epsilon \mathbb{R}$ i $\xi_s=\lambda_s+\theta (\lambda-\lambda_s) \ (|\theta|\leq1)$
for all  $\lambda$ satisfy the condition
  \begin{equation}\label{3.2}
  |\lambda-\lambda_s|<\left|\frac{Q^\prime (\lambda_s)}{Q^{\prime\prime}(\xi_s}\right|
  \end{equation}
  the inequality is true
  \begin{equation}\label{3.3}
  Q(\lambda)>\frac{|\lambda-\lambda_s|}{2}Q^\prime (\lambda_s)
  \end{equation}
 Because
  \begin{eqnarray}\label{3.4}
Q^\prime (\lambda)  =
(a+b)[\cos\lambda(a+b)-kq\cos q\lambda(a+b)]\nonumber\\
Q^{\prime\prime}(\lambda)  =  -(a+b)^2[\sin\lambda(a+b)-kq^2\sin\lambda q(a+b)]
\end{eqnarray}
  then \begin{equation}\label{3.5}|Q^{\prime\prime} (\lambda)|  \leq
(a+b)^2(1+|kq^2|) <(a+b)^2(1+|k|)\end{equation}
  To get a lower estimate for the  $|Q^{\prime} (\lambda_s)| $ we use the  (\ref{3.4}), then we get
  $$\left(Q\prime(w)\right)^2=(a+b)^2\left\{\cos^2w-2kq\cos w\cdot\cos qw+k^2q^2\cos2qw\right\}=(a+b)^2\cdot$$ $$ \cdot\left\{1-\sin^2w+k^2q^2(1-\sin^2qw)-2kq\cos qw \cos w\right\},$$
   where $w=\lambda(a+b)$ and $\sin w=k\sin qw.$ This implies that
 \begin{equation}\label{3.6}  \begin{split}\left(Q\prime(w)\right)^2\geq (a+b)^2\left\{1+k^2q^2-\sin^2w(1+q^2)-2|kq|\sqrt{\left(1-\sin^2w\right)\left(1-\sin^2qw\right)}\right\}\geq \\ \geq (a+b)^2\left\{1+k^2q^2-\sin^2w(1+q^2)-2|kq|\left(1-k^2\sin^2qw\right)\right\}\geq (a+b)^2\left\{1-|kq|^2-\sin^2w(1-|q|^2)\right\}\geq \\ \geq (a+b)^2\left(|q|(1-|k|)\right)(2-|q|-|qk|)>2(a+b)^2|q|(1-|q|)(1-|k|).
   \end{split}
   \end{equation}
 Then
  \begin{equation}\label{3.7}  \begin{split}\left|Q\prime(\lambda_s(0))\right|> \sqrt{2} (a+b)\sqrt{|q|(1-|k|)(1-|q|)}> \\ > (a+b)|q|(1-|q|)(1-|k|)=(a+b)|q|r,
   \end{split}
   \end{equation}
  where
  \begin{equation}\label{3.8}
  r=(1-|q|)(1-|k|)=4\frac{min(a,b)min(1,p^2)}{(a+b)(p^2+1)}<1,
   \end{equation}
  Based on  (\ref{W}) therefore, according to  (\ref{3.7}), (\ref{3.8}) the inequality (\ref{3.2})  is certainly satisfied if
  $$|\lambda-\lambda_s|<\frac{|q|r}{(a+b)(1+|k|)}$$
\begin{Lem1}{}
For all real  $\lambda,$ from the neighborhood
\begin{equation}\label{3.9}
|\lambda-\lambda_s|<\frac{|q|r}{(a+b)(1+|k|)}=R
\end{equation}
  of the zero $\lambda_s(0)$ of the function $\bigtriangleup(0,\lambda)$ (\ref{Det}),  the inequality is valid
  \begin{equation}\label{3.10}
 | \bigtriangleup(0,\lambda)|>\frac{|\lambda-\lambda_s(0)|}{2}|\lambda|(1+p^2)Q^\prime(\lambda_s(0))>\frac{|\lambda-\lambda_s(0)|}{2}|\lambda|
 (1+p^2)(a+b)|q|r,
  \end{equation}
   where $r,\ q$ has form  (\ref{W}),\ (\ref{3.8})
\end{Lem1}

  It follows from the  (\ref{delta}) that
  $$
 |\bigtriangleup(q,\lambda)|>|\bigtriangleup(0,\lambda)|-|\Phi(\lambda)|.
  $$
   We choose  $\lambda\epsilon R$ from the neighborhood (\ref{3.9}) $|\lambda-\lambda_s(0)|<R$ of the zero  $\lambda_s(0) \ (\neq0)$ of the function  $\bigtriangleup(0,\lambda),$ then using  (\ref{phi}) $(\beta=0)$ and (\ref{3.10}) we obtain that
   $$|\bigtriangleup(q,\lambda)|>\frac{|\lambda-\lambda_s(0)|}{2}|\lambda|(1+p^2)Q^\prime(\lambda_s(0))-\delta_1|\lambda|-\delta_2=|\lambda|
   \left(\frac{|\lambda-\lambda_s(0)|}{2}(1+p^2)Q^\prime(\lambda_s(0))-\delta_1-\frac{\delta_2}{|\lambda|}\right), $$
     where numbers  $\delta_s-$ has form (\ref{2.19}). Therefore $|\lambda-\lambda_s|<R$ (\ref{3.9}), then
     $$|\lambda|>|\lambda_s|-R>|\lambda_1|-R>\frac{\pi}{2(a+b)}-\frac{|q|r}{(a+b)(1+|k|)}>\frac{1}{a+b}\left(\frac{\pi}{2}-r\right)>0$$
     based on  {remark 2} , and that mean
     $$|\bigtriangleup(q,\lambda)|>|\lambda|
   \left(\frac{|\lambda-\lambda_s(0)|}{2}(1+p^2)Q^\prime(\lambda_s(0))-\delta_1-\frac{\delta_2(a+b)}{\frac{\pi}{2}-r}\right)$$
     if the first part of this inequality is greater than zero, then
      $$|\lambda-\lambda_s(0)|>\frac{2\delta_1+\frac{4\delta_2(a+b)}{\pi-2r}}{(1+p^2)Q^\prime(\lambda_s(0))}$$
      then for such  $\lambda\epsilon R$ function $|\bigtriangleup(q,\lambda)|$ does not turn to zero. So, if
      \begin{equation}\label{3.11}
      \frac{2\delta_1+\frac{4\delta_2(a+b)}{\pi-2r}}{(1+p^2)Q^\prime(\lambda_s(0))}<|\lambda-\lambda_s(0)|<R,
      \end{equation}
      then  $|\bigtriangleup(q,\lambda)|\neq0$  multiplicity (\ref{3.11})isn`t empty, if
      $$\frac{2\delta_1+\frac{4\delta_2(a+b)}{\pi-2r}}{(1+p^2)Q^\prime(\lambda_s(0))}<R,$$
      and using  (\ref{3.7}) i (\ref{3.9}),  we find that this inequality will certainly be satisfied if
  \begin{equation}\label{3.12}
      2\delta_1+\frac{4\delta_2(a+b)}{\pi-2r}<(1+p^2)\frac{q^2r^2}{1+|k|}
      \end{equation}
      So if the  $\delta_1$ and $\delta_2$ (\ref{2.19}) are such that holds  (\ref{3.12}), then the function  $\bigtriangleup(q,\lambda)$ on the multiplicity  (\ref{3.11}) does not turn to 0. The signs $\bigtriangleup(q,\lambda)$ and $\bigtriangleup(0,\lambda)$ on the left and right sides of multiplicity  (\ref{3.11}) coincide, and given that the signs of the function $\bigtriangleup(0,\lambda)$ on these parts are different, it follows that  $\bigtriangleup(q,\lambda)$ it has at least one root on the multiplicity.
      $$|\lambda-\lambda_s(0)|<\frac{2\delta_1+\frac{4\delta_2(a+b)}{\pi-2r}}{(1+p^2)Q^\prime(\lambda_s(0))}$$
 \begin{Lem1}{}
If numbers  $\delta_1$ and $\delta_2$ (\ref{2.19}) satisfy inequality  (\ref{3.12}), where $p,\ q,\ r$ has form  (\ref{W}) and (\ref{3.8}), then in the surrounding area
\begin{equation}\label{3.13}
|\lambda-\lambda_s(0)|<\frac{2\delta_1+\frac{4\delta_2(a+b)}{\pi-2r}}{(1+p^2)(a+b)|q|r}
  \end{equation}
the zeros  $\lambda_s(0)$ of the function  $\bigtriangleup(0,\lambda)$ (\ref{Det}) contains at least one root  $\lambda_s(q),$ of the perturbed characteristic function  $\bigtriangleup(q,\lambda)$ (\ref{2.19}).
\end{Lem1}
\section{}
  To prove that the characteristic function $\bigtriangleup(q,\lambda)$has no other zeros, except $\lambda_s(q)$  we use Rousche's theorem. Let us denote by  $\gamma_l$ the contour in the  $\mathbb{C},$ formed by the straight lines that connect the points  $\pi\frac{l}{a+b}(1+i), \ \pi\frac{l}{a+b}(-1+i), \ \pi\frac{l}{a+b}(-1-i), \ \pi\frac{l}{a+b}(1-i), (l\epsilon \mathbb{N}). $
  We need a lower estimate for the function  $\bigtriangleup(0,\lambda)$ ) on the contour $\gamma_l$  or, taking into account (\ref{3.1}) a lower estimate for the function  $Q(\lambda)$. For  $\lambda=\alpha+i\beta\ \epsilon\mathbb{C} \ (c=a+b)$ have
    $$Q(\lambda)=\sin(\alpha+i\beta)c-k\sin q(\alpha+i\beta)c=\sin\alpha c \cosh\beta c+i\cos\alpha c\sinh\beta c-k(\sin\alpha qc\cosh\beta qc+i\cos\alpha qc\sinh\beta qc),$$
  then
  $$|Q(\lambda)|^2=\sin^2\alpha c \cosh^2\beta c+k^2\sin^2\alpha qc\cosh^2\beta qc -2k\sin\alpha c\sin\alpha qc\cosh\beta qc\cosh\beta c+\cos^2\alpha c\sinh^2\beta c+$$
  $$+k^2\cos^2\alpha qc\sinh^2\beta qc-2k\cos\alpha c\cos\alpha qc\sinh\beta c\sinh\beta qc=\cosh^2\beta c-\cos^2\alpha c+k^2(\cosh^2\beta qc-\cos^2\alpha qc)-$$
  $$-2k\sin\alpha c\sin\alpha qc\cosh\beta qc\cosh\beta c-2k\cos\alpha c\cos\alpha qc\sinh\beta c\sinh\beta qc\geq(\cosh\beta c-|k|\cosh\beta qc)^2-$$
  $$-(\cos^2\alpha c+k^2\cos^2\alpha qc)(1+|\sinh\beta c||\sinh\beta qc|)\geq(\cosh\beta c-|k|\cosh\beta qc)^2-(1+k^2)(1+|\sinh\beta c||\sinh\beta qc|).$$
    It follows that
$$|Q(\lambda)|\geq(\cosh\beta c-|k|\cosh\beta qc)\sqrt{1-(1+k^2)\frac{(1+|\sinh\beta c||\sinh\beta qc|)}{(\cosh\beta c-|k|\cosh\beta qc)^2}}$$
Hence follows the statement
 \begin{Lem1}{}
At  $\lambda=\alpha+i\beta\epsilon \mathbb{C}$ for function  $\bigtriangleup(0,\lambda)$  (\ref{3.1}) the inequality is true
\begin{equation}\label{4.1}
\begin{split}
|\bigtriangleup(0,\lambda)|>|\lambda|(p+1)\cosh\beta q(a+b)\sqrt{1+k^2}\cdot\\ \cdot\left(1-|\sin\alpha(a+b)\sin\alpha q(a+b)|-
\left(\cos^2\alpha(a+b)+k^2\cos^2\alpha q(a+b)\right)\frac{1+\cosh^2\beta(a+b)}{\cosh^2\beta q(a+b)(1+k^2)}\right)^{1/2}
\end{split}
  \end{equation}
\end{Lem1}
Through $\gamma_l$ we denote the contour in  $\mathbb{C}$ formed by the square with the vertices at the points $\pi\frac{l}{a+b}(1+i), \ \pi\frac{l}{a+b}(-1+i), \ \pi\frac{l}{a+b}(-1-i), \ \pi\frac{l}{a+b}(1-i), (l\epsilon \mathbb{N}). $  On the vertical section (\ref{4.1}) $\lambda=\frac{\pi l}{a+b}(1+\beta i) \ (-1<\beta<1)$ it follow that
$$|\bigtriangleup(0,\lambda)|>\frac{\pi l}{a+b}\sqrt{1+\beta^2}|p+1|\sqrt{1+k^2}\cosh\beta q(a+b)\left(1+\frac{1+\cosh^2\beta(a+b)}{\cosh^2\beta q(a+b)}\right)^{1/2},$$
 and from theorem 3 it follows that for such  $\lambda$ we have
 $$|\Phi(\lambda)|<\cosh\beta a\cosh\beta b(\delta_1|\lambda|+\delta_2),$$
  then at  $l\gg 1$ for $\forall\lambda=\frac{\pi l}{a+b}(1+\beta i) \ \beta\epsilon[-1,1]$ we have
 \begin{equation}\label{4.2}
|\bigtriangleup(0,\lambda)|>|\Phi(\lambda)|
  \end{equation}
It is proved in a similar way that on the sides of the square  $\gamma_l$ at $l\gg 1$  the inequality is true (\ref{4.2}).
\begin{Th}
Suppose that the functions  $q_1(x)$ and $q_2(x)$ in (\ref{Operator}) are such that inequality (\ref{3.12}) holds, where $p, q, r$ are of the form  (\ref{W}) and (\ref{3.8}). Then in each neighborhood  (\ref{3.13})  of the zero  $\lambda_s(0)$ of the characteristic function  $\bigtriangleup(0,\lambda)$ (\ref{Det}) of the unperturbed operator  $L_0$ there is only one zero  $\lambda_s(q)$ of the perturbed characteristic function $\bigtriangleup(q,\lambda)$ (\ref{2.19}) of the operator  $L_q.$
\end{Th}

Therefore, when the potentials are small  $q_1(x)$ and $q_2(x)$ which are expressed only in terms of the parameters of the boundary conditions (\ref{Domain}) each corresponding value of the operator  $L_q$  is located in a small neighborhood of the corresponding value of the unperturbed value of the operator  $L_0.$

\renewcommand{\refname}{\normalsize \rm 

\end{document}